\documentclass[12pt]{article}%
\usepackage{amsmath,amssymb}
\usepackage[applemac]{inputenc}
\usepackage{amsmath,amssymb,fullpage}
\usepackage{color}
\usepackage{amsmath}
\usepackage{amsfonts}
\usepackage{amssymb}
\usepackage{graphicx}
\usepackage{enumerate}
\providecommand{\U}[1]{\protect\rule{.1in}{.1in}}

\newcommand{\fproof}{\hfill $\square$ \bigskip}
\newtheorem{definition}{Definition}[section]

\newtheorem{theorem}[definition]{Theorem}
\newtheorem{problem}[definition]{Problem}
\newtheorem{remark}[definition]{ \it Remark}

\newtheorem{proposition}[definition]{Proposition}
\newtheorem{lemma}[definition]{Lemma}

\numberwithin{equation}{section}

\def\1B{\text{1\!\!I}}

\begin{document}

\title{Singular control of stochastic Volterra integral equations }
\author{Nacira Agram$^{1}$, Saloua Labed$^{2}$, Bernt \O ksendal$^{3}$ and Samia
Yakhlef$^{2}$ }
\date{14 April 2021}
\maketitle

\begin{abstract}
{\small \noindent This paper deals with optimal combined singular and regular
controls for stochastic Volterra integral equations, where the solution
$X^{u,\xi}(t)=X(t)$ is given by
\[%
\begin{array}
[c]{cc}%
X(t) & =\phi(t)+{\textstyle\int_{0}^{t}}b\left(  t,s,X(s),u(s)\right)  ds+%
{\textstyle\int_{0}^{t}}
\sigma\left(  t,s,X(s),u(s)\right)  dB(s)\\
& +%
{\textstyle\int_{0}^{t}}
h\left(  t,s\right)  d\xi(s).
\end{array}
\]
Here $dB(s)$ denotes the Brownian motion It\^o type differential and $\xi$ denotes the singular control (singular in time $t$ with respect to Lebesgue measure) and $u$ denotes the regular control (absolutely continuous with respect to Lebesgue measure). \newline Such systems may for example be used to model 
harvesting of populations with memory, where $X(t)$ represents the population
density at time $t$, and the singular control process $\xi$ represents the
harvesting effort rate. The total income from the harvesting is represented
by
\[%
\begin{array}
[c]{ll}%
J(u,\xi) & =\mathbb{E[}{%
{\textstyle\int_{0}^{T}}
}\text{ }f_{0}(t,X(t),u(t))dt+{%
{\textstyle\int_{0}^{T}}
}f_{1}(t,X(t))d\xi(t)+g(X(T))],
\end{array}
\]
for given functions $f_{0},f_{1}$ and $g$, where $T>0$ is a constant denoting
the terminal time of the harvesting. Note that it is important to allow the controls to be singular, because in some cases the optimal controls are of this type. \newline}

{\small \noindent Using Hida-Malliavin calculus, we prove sufficient
conditions and necessary conditions of optimality of controls. As a
consequence, we obtain a new type of backward stochastic Volterra integral
equations with singular drift. \newline}

{\small \noindent Finally, to illustrate our results, we apply them to discuss
optimal harvesting problems with possibly density dependent prices. }

\end{abstract}

\footnotetext[1]{{\small Department of Mathematics, Linnaeus University (LNU), V\" axj\" o, Sweden. Email:\texttt{ nacira.agram@lnu.se.}}}

\footnotetext[2]{{\small University Mohamed Khider of Biskra, Algeria. Email:
\texttt{s.labed@univ-biskra.dz, s.yakhlef@univ-biskra.dz}}}

\footnotetext[3]{{\small Department of Mathematics, University of Oslo, P.O.
Box 1053 Blindern, N--0316 Oslo, Norway. Email:\texttt{ oksendal@math.uio.no.}%
}}

\paragraph{{\protect\small MSC(2010):}}

{\small 60H05, 60H20, 60J75, 93E20, 91G80,91B70.}

\paragraph{{\protect\small Keywords:}}

{\small Stochastic maximum principle; stochastic Volterra integral equation; singular control; backward stochastic Volterra integral equation; Hida-Malliavin
calculus.}

\section{Introduction}

As a motivating example, consider the population of a certain type of
fish in a lake, where the density $X(t)$ at time $t$ can be modelled as the solution of the following stochastic
Volterra integral equation (SVIE):
\[
X(t)=x_{0}+%
{\textstyle\int_{0}^{t}}
b_{0}(t,s)X(s)ds+{%
{\textstyle\int_{0}^{t}}
}\sigma_{0}(s)X(s)dB(s)-%
{\textstyle\int_{0}^{t}}
\gamma_{0}(t,s)d\xi(s),
\]
where  the coefficients
$b_{0},\sigma_{0}$ {and} $\gamma_{0} $ are bounded deterministic functions,
and \newline $B(t)=\{B(t,\omega)\}_{t\geq0,\omega \in \Omega}$ {\small is a Brownian motion defined on a complete
probability space} $(\Omega,\mathcal{F},P).$ {\small We associate to this
space a natural filtration} $\mathbb{F}=\{\mathcal{F}_{t}\}_{t\geq0}$
generated by $B(t)$, assumed to satisfy the usual conditions. The 
process $\xi(t)$ is our control process. It is an $\mathbb{F}$- adapted,
nondecreasing left-continuous process representing the harvesting effort. It is called singular, because as a function of time $t$ it may be singular with respect to Lebesgue measure. The
constant $\gamma_{0}>0$ is the harvesting efficiency coefficient. It turns out
that in some cases the optimal process $\xi(t)$ can be represented as the
local time of the solution $X(t)$ at some threshold curve. In this case $\xi(t)$ is increasing only on a set of Lebesgue measure $0$.

\noindent Volterra equations are commonly used in population growth models,
especially when age dependence plays a role. See e.g. Gripenberg \textit{et
al} \cite{GLS}. Moreover, they are important examples of equations with memory.

\noindent We assume that the total expected utility from the harvesting is
represented by
\[
J(\xi)=\mathbb{E[}\theta X(T)+%
{\textstyle\int_{0}^{T}}
\log(X(t))d\xi(t)],
\]
where $\mathbb{E}$ denotes the expectation with respect to $P$. The problem is
then to maximise $J(\xi)$ over all admissible singular controls $\xi$.
We will return to this example in Section 4.\newline

\noindent Control problems for singular Volterra integral equations have been
studied by Lin and Yong \cite{LY} in the deterministic case. In this paper we
study singular control of SVIEs and we present a different approach based on a stochastic
version of the Pontryagin maximum principle.

\noindent Stochastic control for Volterra integral equations has been studied
by Yong \cite{Yong1} and subsequently by Agram \textit{el al }\cite{AO},
\cite{AOY} who used the white noise calculus to obtain both sufficient and
necessary conditions of optimality. In the latter, smoothness of coefficients
is required.

\noindent The adjoint processes of our maximum principle satisfy a backward
stochastic integral equation of Volterra type and with a singular term coming
from the control. In our example one may consider the optimal singular term as
the local time of the state process that is keeping it above/below a certain
threshold curve. Hence in some cases we can associate this type of equations
with reflected backward stochastic Volterra integral equations.

\noindent Partial result for existence and uniqueness of backward stochastic
Volterra integral equation (BSVIE) in a continuous case can be found in Yong \cite{Yong1},
\cite{Yong 2}, and for a discontinuous case, we refer for example to Agram \textit{el al }\cite{AOY2}, \cite{A} where there
are also some applications.
\newline 

\noindent The paper is organised as follows: In the next section we give some
preliminaries about the generalised Malliavin calculus, called
Hida-Malliavin calculus, in the white noise space of Hida of stochastic
distributions. Section 3 is addressed to the study of the stochastic maximum
principle where both sufficient and necessary conditions of optimality are
proved. Finally, in Section 4 we apply the results obtained in section 3 to
discuss optimal harvesting problems with possibly density dependent prices.

\section{Hida - Malliavin calculus}

{\small \noindent Let }$\mathbb{G}=\{\mathcal{G}_{t}\}_{t\geq0}$ be a
subfiltration of $\mathbb{F}$, {\small in the sense that} $\mathcal{G}%
_{t}\subseteq\mathcal{F}_{t},$ {\small for all} $t\geq0.$ {\small The given
set }$U\subset\mathbb{R}$ {\small is assumed to be convex. The set of
admissible controls, i.e. the strategies available to the controller, is given
by a subset} $\mathcal{A}$ {\small of the c\`adl\`ag, }$U${\small -valued and
}$\mathbb{G}${\small -adapted processes. Let }$\mathcal{K}$ {\small be the set
of all }$\mathbb{G}${\small -adapted processes }$\xi(t)$ {\small that are
nondecreasing and left continuous with respect to }$t${\small .}

{\small \noindent Next we present some preliminaries about the extension of
the Malliavin calculus into the stochastic distribution space of Hida, for
more details, we refer the reader to Aase \textit{et al} \cite{AaOPU}, Di
Nunno \textit{et al} \cite{DOP}.}

{\small \noindent The classical Malliavin derivative is only defined on a
subspace} $\mathbb{D}_{1,2}$ {\small of} $\mathbb{L}^{2}(P)${\small . However,
there are many important random variables in $\mathbb{L}^{2}(P)$ that do not
belong to $\mathbb{D}_{1,2}$. For example, this is the case for the solutions
of a backward stochastic differential equations or more generally the BSVIE.
This is why the Malliavin derivative was extended to an operator defined on
the whole of} $\mathbb{L}^{2}(P)$ {\small and with values in the Hida space}
$(\mathcal{S})^{\ast}$ {\small of stochastic distributions. It was proved by
Aase \textit{et al} \cite{AaOPU} that one can extend the Malliavin derivative
operator }$D_{t}$ {\small from }$\mathbb{D}_{1,2}$ {\small to all of}
$\mathbb{L}^{2}(\mathcal{F}_{T},P)$ {\small in such a way that, also denoting
the extended operator by }$D_{t}${\small , for all random variable }%
$F\in\mathbb{L}^{2}(\mathcal{F}_{T},P)${\small , we have}
\begin{equation}
D_{t}F\in(\mathcal{S})^{\ast}\text{ {\small and} }(t,\omega)\mapsto
\mathbb{E}[D_{t}F|\mathcal{F}_{t}]\text{ {\small belongs to} }\mathbb{L}%
^{2}(\lambda\times P), \label{eq2.10a}%
\end{equation}
{\small \noindent where }$\lambda$ {\small is Lebesgue measure on }$[0,T].$ We
now give a short introduction to Malliavin calculus and its extension to Hida-Malliavin calculus in the white noise setting:

\begin{definition}
\label{direct.derivative}

\begin{description}
\item[(i)] {\small Let }$F\in\mathbb{L}^{2}(P)$ {\small and let }$\gamma
\in\mathbb{L}^{2}(\mathbb{R})$ {\small be deterministic. Then the directional
derivative of }$F$ {\small in} $(\mathcal{S})^{\ast}$ {\small (respectively,
in} $\mathbb{L}^{2}(P)${\small ) in the direction }$\gamma$ {\small is defined
by }%
\begin{equation}
D_{\gamma}F(\omega)=\lim_{\varepsilon\rightarrow0}\frac{1}{\varepsilon
}\big[F(\omega+\varepsilon\gamma)-F(\omega)\big] \label{6.1}%
\end{equation}
{\small whenever the limit exists in }$(\mathcal{S})^{\ast}$
{\small (respectively, in }$\mathbb{L}^{2}(P)${\small ).}

\item[(ii)] Suppose there exists a function $\psi:\mathbb{R}\mapsto
(\mathcal{S})^{\ast}$ (respectively, $\psi:\mathbb{R}\mapsto\mathbb{L}^{2}%
(P)$) such that
\begin{equation}%
\begin{split}
&
{\textstyle\int_{{\mathbb{R}}}}
\psi(t)\gamma(t)dt\quad\text{ exists in }(\mathcal{S})^{\ast}\text{
(respectively, in }\mathbb{L}^{2}(P))\text{ and}\\
&  D_{\gamma}F=%
{\textstyle\int_{\mathbb{R}}}
\psi(t)\gamma(t)dt,\quad\text{ for all }\gamma\in\mathbb{L}^{2}(\mathbb{R}).
\end{split}
\label{6.2}%
\end{equation}
Then we say that $F$ is \emph{Hida-Malliavin differentiable} in $(\mathcal{S}%
)^{\ast}$ (respectively, in $\mathbb{L}^{2}(P)$) and we write
\[
\psi(t)=D_{t}F,\quad t\in\mathbb{R}.
\]
We call $D_{t}F$ \label{simb-053}the \emph{Hida-Malliavin derivative at $t$ in
$(\mathcal{S})^{\ast}$ (respectively, in }$\mathbb{L}$\emph{$^{2}(P)$)} or the
\emph{stochastic gradient} of $F$ at $t$.
\end{description}
\end{definition}

\noindent Let $F_{1},...,F_{m}\in\mathbb{L}^{2}(P)$ be Hida-Malliavin
differentiable in $\mathbb{L}^{2}(P)$. Suppose that $\varphi\in C^{1}%
({\mathbb{R}}^{m})$, $D_{t}F_{i}\in\mathbb{L}^{2}(P)$, for all $t\in
{\mathbb{R}}$, and $\frac{\partial\varphi}{\partial x_{i}}(F)D_{\cdot}F_{i}%
\in\mathbb{L}^{2}(\lambda\times P)$ for $i=1,...,m$, where $F=(F_{1}%
,...,F_{m})$. Then $\varphi(F)$ is Hida-Malliavin differentiable and
\begin{equation}
D_{t}\varphi(F)=%
{\textstyle\sum_{i=1}^{m}}
\tfrac{\partial\varphi}{\partial x_{i}}(F)D_{t}F_{i}. \label{7.4bis}%
\end{equation}

\noindent We have the following \emph{generalized duality formula,} for the
Brownian motion:

\begin{proposition}
Fix $s\in\lbrack0,T]$. If $t\mapsto\varphi(t,s,\omega)\in\mathbb{L}%
^{2}(\lambda\times P)$ is $\mathbb{F}$-adapted with \newline$\mathbb{E}[%
{\textstyle\int_{0}^{T}}
\varphi^{2}(t,s)dt]<\infty$ and $F\in\mathbb{L}^{2}(\mathcal{F}_{T},P)$, then
we have
\begin{equation}
\mathbb{E}[F%
{\textstyle\int_{0}^{T}}
\varphi(t,s)dB(t)]=\mathbb{E}[%
{\textstyle\int_{0}^{T}}
\mathbb{E}[D_{t}F|\mathcal{F}_{t}]\varphi(t,s)dt]. \label{geduB}%
\end{equation}

\end{proposition}

We will need the following:

\begin{lemma}
Let $t,s,\omega\mapsto G(t,s,\omega)\in\mathbb{L}^{2}(\lambda\times
\lambda\times P)$ and $t,\omega\mapsto p(t)\in\mathbb{L}^{2}(\lambda\times
P),$ then the followings hold:

\begin{enumerate}
\item The Fubini theorem combined with a change of variables gives
\begin{equation}%
{\textstyle\int_{0}^{T}}
p(t)(%
{\textstyle\int_{0}^{t}}
G(t,s)ds)dt={\textstyle\int_{0}^{T}} ({\textstyle\int_{t}^{T}}
p(s)G(s,t)ds)dt, \label{1}%
\end{equation}
and
\begin{equation}
{\textstyle\int_{0}^{T}} p(t)({\textstyle\int_{0}^{t}} G(t,s)ds)d\xi(t)
={\textstyle\int_{0}^{T}} ({\textstyle\int_{t}^{T}} p(s)G(s,t)ds)d\xi(t).
\label{2}%
\end{equation}

\item[2.] The generalized duality formula $\left(  \ref{geduB}\right)  $
together with the Fubini theorem, yields
\begin{equation}
\mathbb{E}[%
{\textstyle\int_{0}^{T}}
p(t)(%
{\textstyle\int_{0}^{t}}
G(t,s)dB(s))dt]=\mathbb{E}[
{\textstyle\int_{0}^{T}}
{\textstyle\int_{t}^{T}}
\mathbb{E}[D_{t}p(s)|\mathcal{F}_{t}]G(s,t)dsdt]. \label{3}%
\end{equation}

\end{enumerate}
\end{lemma}

\section{Stochastic maximum principles}

\noindent{\small In this section, we study stochastic maximum principles of
stochastic Volterra integral systems under partial information, i.e., the
information available to the controller is given by a sub-filtration
$\mathbb{G}.$ Suppose that the state of our system $X^{u,\xi}(t)=X(t)$
satisfies the following SVIE%
\begin{equation}%
\begin{array}
[c]{cc}%
X(t) & =\phi(t)+%
{\textstyle\int_{0}^{t}}
b\left(  t,s,X(s),u(s)\right)  ds+%
{\textstyle\int_{0}^{t}}
\sigma\left(  t,s,X(s),u(s)\right)  dB(s)\\
& +%
{\textstyle\int_{0}^{t}}
h\left(  t,s\right)  d\xi(s),\quad t\in\lbrack0,T],
\end{array}
\label{sde}%
\end{equation}
where $b(t,s,x,u)=b(t,s,x,u,\omega):\left[  0,T\right]  ^{2}\times
\mathbb{R}\times U\times\Omega\rightarrow\mathbb{R}$, $\sigma(t,s,x,u)=\sigma
(t,s,x,u,\omega):\left[  0,T\right]  ^{2}\times\mathbb{R}\times U\times
\Omega\rightarrow%
\mathbb{R}
.$\newline\newline The \emph{performance functional} has the form%
\begin{equation}%
\begin{array}
[c]{ll}%
J(u,\xi) & =\mathbb{E[}{%
{\textstyle\int_{0}^{T}}
}\text{ }f_{0}(t,X(t),u(t))dt+{%
{\textstyle\int_{0}^{T}}
}f_{1}(t,X(t))d\xi(t)+g(X(T))],
\end{array}
\label{perf}%
\end{equation}
with given functions $f_{0}(t,x,u)=f_{0}(t,x,u,\omega):\left[  0,T\right]
\times\mathbb{R}\times U\times\Omega\rightarrow\mathbb{R}$}$,$ {\small $f_{1}%
(t,x)=f_{1}(t,x,\omega):\left[  0,T\right]  \times\mathbb{R}\times
\Omega\rightarrow\mathbb{R}$ and $g(x)=g(x,\omega):\mathbb{R\times}%
\Omega\rightarrow\mathbb{R}$. Let $\mathcal{A},\mathcal{K}$ denote the family of admissible controls $u,\xi$, respectively.
We let $\mathcal{A}$ be the set of all adapted {\small of the c\`adl\`ag, } processes $u(t,\omega) \in L^2(dt \times dP)$ and $\mathcal{K}$ consists of all adapted nondecreasing processes $\xi(t)$ with $\xi(0) = 0$.
We study the following problem:
\begin{problem}
Find a control pair
$(\hat{u},\hat{\xi}) \in \mathcal{A} \times \mathcal{K}$ such that
\begin{align}
J(\hat{u},\hat{\xi}) = \sup_{(u,\xi) \in \mathcal{A} \times \mathcal{K}} J(u,\xi).
\end{align} \label{pr3.3} 
\end{problem}

We impose the following assumptions on the coefficients:\newline}
{\small The processes $b(t,s,x,u),\sigma(t,s,x,u),f_{0}(s,x,u),$
$f_{1}(t,x,\xi)$\ and $h(t,s)$ are $\mathbb{F}$-adapted with respect to $s$
for all $s\leq t,$ and twice continuously differentiable ($C^{2}$) with
respect to $t$, $x,$ and continuously differentiable ($C^{1}$) with respect to
$u$ for each $s$. The driver $g$ is assumed to be $\mathcal{F}_{T}$-measurable
and ($C^{1}$) in $x$. Moreover, all the partial derivatives are supposed to be
bounded. }

{\small \noindent Note that the performance functional (\ref{perf}) is not of
Volterra type. }

\subsection{The Hamiltonian and the adjoint equations}

{\small \noindent Define the \emph{Hamiltonian functional} associated to our
control problem (\ref{sde}) and (\ref{perf}), as
\begin{equation}%
\begin{array}
[c]{l}%
\mathbb{H}(t,x,u,p,q)(dt,d\xi(t))\\
:=\big[H_{0}(t,x,u,p,q)+H_{1}(t,x,u,p)\big]dt+\big[\overline{H_{0}%
}(t,x,p)+\overline{H_{1}}(t,p)\big]d\xi(t),
\end{array}
\label{eq3.3}%
\end{equation}
where
\begin{align*}
H_{0}  &  :[0,T]\times\mathbb{R}\times U\times\mathbb{R}\times\mathbb{R}%
\rightarrow\mathbb{R},\\
H_{1}  &  :[0,T]\times\mathbb{R}\times U\times\mathbb{R}^{[0,T]}\rightarrow\mathbb{R},\\
\overline{H_{0}}  &  :[0,T]\times\mathbb{R}\times\mathbb{R}\rightarrow\mathbb{R},\\
\overline{H_{1}}  &  :[0,T]\times\mathbb{R}^{[0,T]}\rightarrow\mathbb{R},
\end{align*}
are defined as follows
\begin{align*}
H_{0}(t,x,u,p,q)  &  :=f_{0}(t,x,u)+p(t)b(t,t,x,u)+q(t,t)\sigma(t,t,x,u),\\
H_{1}(t,x,u,p)  &  :=%
{\textstyle\int_{t}^{T}}
p(s)\tfrac{\partial b}{\partial s}(s,t,x,u)ds+%
{\textstyle\int_{t}^{T}}
\mathbb{E}[D_{t}p(s)|\mathcal{F}_{t}]\tfrac{\partial\sigma}{\partial
s}(s,t,x,u)ds,\\
\overline{H_{0}}(t,x,p)  &  :=f_{1}(t,x)+p(t)h(t,t),\\
\overline{H_{1}}(t,p)  &  :=%
{\textstyle\int_{t}^{T}}
p(s)\tfrac{\partial h}{\partial s}(s,t)ds.
\end{align*}
For convenience, we will use the following notation from now on:
\begin{equation}
\mathcal{H}(t,x,u,p,q)=H_{0}(t,x,u,p,q)+H_{1}(t,x,u,p), \label{h}%
\end{equation}%
\begin{equation}
\overline{\mathcal{H}}(t,x,p)=\overline{H_{0}}(t,x,p)+\overline{H_{1}}(t,p).
\label{h1}%
\end{equation}
The BSVIE for the adjoint processes $p(t),q(t,s)$ is defined by%
\begin{equation}
p(t)=\tfrac{\partial g}{\partial x}(X(T))+%
{\textstyle\int_{t}^{T}}
\tfrac{\partial\mathcal{H}}{\partial x}(s)ds+%
{\textstyle\int_{t}^{T}}
\tfrac{\partial\overline{\mathcal{H}}}{\partial x}(s)d\xi(s)-%
{\textstyle\int_{t}^{T}}
q(t,s)dB(s), \label{p}%
\end{equation}
where we have used the simplified notation
\begin{align*}
\tfrac{\partial\mathcal{H}}{\partial x}(t)  &  =\tfrac{\partial\mathcal{H}%
}{\partial x}(t,X(t),u(t),p(t),q(t,t)),\\
\tfrac{\partial\overline{\mathcal{H}}}{\partial x}(t)  &  =\tfrac
{\partial\overline{\mathcal{H}}}{\partial x}(t,X(t),p(t)).
\end{align*}
Note that from equation (\ref{sde}), we get the following equivalent
formulation, for each $(t,s)\in\lbrack0,T]^{2},$
\begin{equation}%
\begin{array}
[c]{ll}%
dX(t) & =\phi^{\prime}(t)dt+b\left(  t,t,X(t),u(t)\right)  dt+(%
{\textstyle\int_{0}^{t}}
\tfrac{\partial b}{\partial t}\left(  t,s,X(s),u(s)\right)  ds)dt+\sigma
\left(  t,t,X(t),u(t)\right)  dB(t)\\
& +(%
{\textstyle\int_{0}^{t}}
\tfrac{\partial\sigma}{\partial t}\left(  t,s,X(s),u(s)\right)
dB(s))dt+h\left(  t,t\right)  d\xi(t)+(%
{\textstyle\int_{0}^{t}}
\tfrac{\partial h}{\partial t}(t,s)d\xi(s))dt.
\end{array}
\label{eq3.6}%
\end{equation}
We assume that the map $t\mapsto q(t,s)$ is ($C^{1}$) for all $s,\omega$ and
moreover,
\[
\mathbb{E}[  {\textstyle\int_{0}^{T}}{\textstyle\int_{0}^{T}}\left(  \tfrac{\partial
q(t,s)}{\partial t}\right)  ^{2}dsdt]  <\infty,
\]
under which we can write the following differential form of equation
(\ref{p}):
\begin{equation}
\left\{
\begin{array}
[c]{l}%
dp(t)=-[\tfrac{\partial\mathcal{H}}{\partial x}(t)dt+\tfrac{\partial
\overline{\mathcal{H}}}{\partial x}(t)d\xi(t)+%
{\textstyle\int_{t}^{T}}
\tfrac{\partial q}{\partial t}(t,s)dB(s)dt]+q(t,t)dB(t),\\
p(T)=\tfrac{\partial g}{\partial x}(X(T)).
\end{array}
\right.  \label{eq3.8}%
\end{equation}
}

\subsection{A sufficient maximum principle}

{\small \noindent We will see under which conditions the couple $(u,\xi)$ is
optimal, i.e. we will prove a sufficient version of the maximum principle
approach (a verification theorem). }

\begin{theorem}
[Sufficient maximum principle]{\small Let $(\widehat{u},\widehat{\xi}%
)\in\mathcal{A\times}\mathcal{K},$ with corresponding solutions }%
$\widehat{{\small X}}${\small $(t),$ $\left(  \widehat{p}(t),\widehat
{q}(t,s)\right)  $ of (\ref{sde}) and (\ref{p}) respectively. Assume that the
functions $x\mapsto g(x)$ and $(x,u,\xi)\mapsto\mathbb{H}(t,x,u,\widehat
{p},\widehat{q})(dt,d\xi(t))$ are concave. Moreover, we impose the following
optimal conditions for each control: }

\begin{itemize}
\item (Maximum condition for $u$)
\begin{equation}
\underset{u\in U}{\sup}\text{ }\mathbb{E[H}(t,\widehat{{\small X}%
}(t),u,\widehat{p}(t),\widehat{q}(t,t))|{\mathcal{G}}_{t}]=\mathbb{E[H}%
(t,\widehat{{\small X}}(t),\widehat{u}(t),\widehat{p}(t),\widehat
{q}(t,t))|{\mathcal{G}}_{t}],\text{ for a.a. }t,\text{ }P\text{-a.s.}
\label{mc-1}%
\end{equation}
where we are using the notation
\begin{align*}
\mathbb{E[H}(t,\widehat{{\small X}}(t),u,\widehat{p}(t),\widehat
{q}(t,t))|{\mathcal{G}}_{t}]: & = \mathbb{E}[\mathcal{H}(t,\widehat
{{\small X}}(t),u,\widehat{p}(t),\widehat{q}(t,t))|{\mathcal{G}}%
_{t}]dt\nonumber\\
&  +\mathbb{E}[\overline{\mathcal{H}}(t,\widehat{{\small X}}(t),\widehat
{p}(t))|{\mathcal{G}}_{t}]d\xi(t),\text{ for a.a. }t,\text{ }P\text{-a.s.}
\end{align*}

\item (Maximum condition for $\xi$) \newline For all $\xi\in\mathcal{K}$ we
have, in the sense of inequality between random measures,
\begin{align}
&  \mathbb{E[H}(t,\widehat{{\small X}}(t),u,\widehat{p}(t),\widehat
{q}(t,t))(dt,d\xi(t))|{\mathcal{G}}_{t}]\nonumber\\
&  \leq\mathbb{E[H}(t,\widehat{{\small X}}(t),u,\widehat{p}(t),\widehat
{q}(t,t))(dt,d\widehat{\xi}(t))|{\mathcal{G}}_{t}],\text{ for each }t,\text{
}P\text{-a.s.} \label{eq3.11}%
\end{align}

\end{itemize}

Then {\small $(\widehat{u},\widehat{\xi})$} is an optimal pair.
\end{theorem}

{\small \noindent{Proof.} \quad Choose $u\in\mathcal{A}$ and $\xi
\in\mathcal{K},$ we want to prove that $J(u,\xi)-J(\widehat{u},\widehat{\xi
})\leq0$. We set%
\[
J(u,\xi)-J(\widehat{u},\widehat{\xi})=J(u,\xi)-J(u,\widehat{\xi}%
)+J(u,\widehat{\xi})-J(\widehat{u},\widehat{\xi}).
\]
Since we have one regular control and one singular, we will solve the problem
by separating them, as follows: }

{\small \noindent First, we prove that $\xi$ is optimal i.e., for all fixed
$u\in U$, $J(u,\xi)-J(u,\widehat{\xi})\leq0.$ Then, we plug the optimal
}$\widehat{\xi}$ {\small into the second part and we prove it for $u$, i.e.,
$J(u,\widehat{\xi})-J(\widehat{u},\widehat{\xi})\leq0.$ However, the case of
regular controls $u$ has been proved in Theorem 4.3 by Agram \textit{et al}
\cite{AOY2} . It rests to prove only the inequality for singular controls
$\xi$. }

{\small \noindent From definition (\ref{perf}), we have
\begin{equation}
J(u,\xi)-J(u,\widehat{\xi})=A_{1}+A_{2}+A_{3}, \label{eq2.8}%
\end{equation}
where we have used hereafter the shorthand notations \newline%
\[
A_{1}=\mathbb{E[}%
{\textstyle\int_{0}^{T}}
\widetilde{f_{0}}\left(  t\right)  dt],\quad A_{2}=\mathbb{E[}%
{\textstyle\int_{0}^{T}}
f_{1}\left(  t\right)  d\xi(t)-%
{\textstyle\int_{0}^{T}}
\widehat{f_{1}}\left(  t\right)  d\widehat{\xi}(t)],\quad A_{3}=\mathbb{E[}%
\widetilde{g}(T)],
\]
with }$\widetilde{{\small f_{0}}}${\small $\left(  t\right)  =f_{0}%
(t)-\widehat{f_{0}}(t)$, }$\widetilde{{\small g}}$%
{\small $(T)=g(X(T))-g(\widehat{X}(T))$, and similarly for $b(t,t)=b\left(
t,t,X(t),u(t)\right)  $, and the other coefficients. By definition $(\ref{h}%
)$, we get
\begin{equation}
A_{1}=\mathbb{E}[%
{\textstyle\int_{0}^{T}}
\{\widetilde{H_{0}}(t)-\widehat{p}(t)\widetilde{b}(t,t)-\widehat
{q}(t,t)\widetilde{\sigma}(t,t)\}dt]. \label{i1}%
\end{equation}
Concavity of $g$ together with the terminal value of the BSVIE $\left(
\ref{p}\right)  $, we obtain%
\[%
\begin{array}
[c]{lll}%
A_{3} & \leq\mathbb{E}[\tfrac{\partial\widehat{g}}{\partial x}(T)\widetilde
{X}(T)] & =\mathbb{E}[\widehat{p}(T)\widetilde{X}(T)].
\end{array}
\]
Applying the integration by parts formula to the product }$\widehat{p}%
${\small $(t)\widetilde{X}(t)$, we get
\begin{align}
A_{3}  &  \leq\mathbb{E}[\widehat{p}(T)\widetilde{X}(T)]\nonumber\\
&  =\mathbb{E[}%
{\textstyle\int_{0}^{T}}
\widehat{p}(t)\{\widetilde{b}(t,t)+%
{\textstyle\int_{0}^{t}}
\tfrac{\partial\widetilde{b}}{\partial t}(t,s)ds+%
{\textstyle\int_{0}^{t}}
\tfrac{\partial\widetilde{\sigma}}{\partial t}\left(  t,s\right)  dB(s)+%
{\textstyle\int_{0}^{t}}
\tfrac{\partial h}{\partial t}(t,s)d\widetilde{\xi}(s)\}dt\nonumber\\
&  +%
{\textstyle\int_{0}^{T}}
\widehat{p}(t)\widetilde{\sigma}(t,t)dB(t)+%
{\textstyle\int_{0}^{T}}
\widehat{p}(t)h(t,t)d\tilde{\xi}(t)-%
{\textstyle\int_{0}^{T}}
\widetilde{X}(t)\tfrac{\partial\widehat{\mathcal{H}}}{\partial x}(t)dt-%
{\textstyle\int_{0}^{T}}
\widetilde{X}(t)\tfrac{\partial\widehat{\overline{\mathcal{H}}}}{\partial
x}(t)d\widehat{\xi}(t)\nonumber\\
&  -%
{\textstyle\int_{0}^{T}}
\widetilde{X}(t)(%
{\textstyle\int_{t}^{T}}
\tfrac{\partial\widehat{q}}{\partial t}(t,s)dB(s))dt+%
{\textstyle\int_{0}^{T}}
\widetilde{X}(t)\widehat{q}(t,t)dB(t)+%
{\textstyle\int_{0}^{T}}
\widehat{q}(t,t)\widetilde{\sigma}(t,t)dt]. \label{I2}%
\end{align}
It follows from formulas (\ref{1})-(\ref{3}), that }

{\small
\[%
\begin{array}
[c]{ll}%
\mathbb{E[}{\int_{0}^{T}}\widehat{p}(t)({\int_{0}^{t}}\tfrac{\partial
\widetilde{b}}{\partial t}(t,s)ds)dt] & =\mathbb{E[}{\int_{0}^{T}}({\int
_{t}^{T}}\widehat{p}(s)\tfrac{\partial\widetilde{b}}{\partial s}(s,t)ds)dt],\\
\mathbb{E[}{\int_{0}^{T}}\widehat{p}(t)({\int_{0}^{t}}\tfrac{\partial
h}{\partial t}(t,s)d\widetilde{\xi}(s))dt] & =\mathbb{E[}{\int_{0}^{T}}%
({\int_{t}^{T}}\widehat{p}(s)\tfrac{\partial h}{\partial s}(s,t)ds)d\widetilde
{\xi}(t)],\\
\mathbb{E[}{\int_{0}^{T}}\widehat{p}(t)({\int_{0}^{t}}\tfrac{\partial
\widetilde{\sigma}}{\partial t}(t,s)dB(s))dt] & =\mathbb{E[}{\int_{0}^{T}%
\int_{t}^{T}\mathbb{E}(D_{t}\widehat{p}(s)|\mathcal{F}_{t})\tfrac
{\partial\widetilde{\sigma}}{\partial s}(s,t,x,u)}dsdt].
\end{array}
\]
Substituting the above into $\left(  \ref{eq2.8}\right)  ,$ we obtain
\begin{align*}
J(\widehat{u},\xi)-J(\widehat{u},\widehat{\xi})  &  \leq\mathbb{E}[{\textstyle\int_{0}^{T}}(\widetilde{H}_{0}(t)+\widetilde{H}_{1}(t))dt+{\textstyle\int_{0}^{T}}
f_{1}(t)d\xi(t)-{\textstyle\int_{0}^{T}}\widehat{f}_{1}(t)d\widehat{\xi}(t)+{\textstyle\int_{0}^{T}}\widehat{p}(t)h(t,t)d\widetilde{\xi}(t)\\
&  +{\textstyle\int_{0}^{T}}({\textstyle\int_{0}^{T}}\widehat{p}(s)\tfrac{\partial h}{\partial
s}(s,t)ds)d\widetilde{\xi}(t)-{%
{\textstyle\int_{0}^{T}}
}\widetilde{X}(t)\tfrac{\partial\widehat{\mathcal{H}}}{\partial x}(t)dt-{%
{\textstyle\int_{0}^{T}}
}\widetilde{X}(t)\tfrac{\partial\widehat{\overline{\mathcal{H}}}}{\partial
x}(t)d\widehat{\xi}(t)]\\
&  =\mathbb{E}[{\textstyle\int_{0}^{T}}(\mathcal{H}(t)-\widehat{\mathcal{H}%
}(t))dt+(\overline{\mathcal{H}}(t)d\xi(t)-\widehat{\overline{\mathcal{H}}%
}(t)d\widehat{\xi}(t))-{\textstyle\int_{0}^{T}}\widetilde{X}(t)\tfrac{\partial
\widehat{\mathcal{H}}}{\partial x}(t)dt\\
& -{\textstyle\int_{0}^{T}}\widetilde{X}%
(t)\tfrac{\partial\widehat{\overline{\mathcal{H}}}}{\partial x}(t)d\widehat
{\xi}(t)].
\end{align*}
Using the concavity of }$\mathcal{H}$ and $\overline{\mathcal{H}}$ with
respect to $x$ and $\xi,$ we have%
\begin{align*}
&  J({\small \widehat{u}},\xi)-{\small J(\widehat{u},\widehat{\xi})}\\
&  \leq\mathbb{E}[{{{\textstyle\int_{0}^{T}}\{}}\widetilde{X}(t)\tfrac{\partial
\widehat{\mathcal{H}}}{\partial x}(t)-\widetilde{X}(t)\tfrac{\partial
\widehat{\mathcal{H}}}{\partial x}(t)\}dt+{\textstyle\int_{0}^{T}}\widetilde
{X}(t)\tfrac{\partial\widehat{\overline{\mathcal{H}}}}{\partial x}%
(t)d{\small \widehat{\xi}}(t)\\
&  -{\textstyle\int_{0}^{T}}\widetilde{X}(t)\tfrac{\partial\widehat{\overline
{\mathcal{H}}}}{\partial x}(t)d{\small \widehat{\xi}}(t)+{\textstyle\int_{0}^{T}}\widehat{\overline{\mathcal{H}}}(t)\{d\xi(t)-d{\small \widehat{\xi}}(t)\}]\\
&  =\mathbb{E}[{\textstyle\int_{0}^{T}}\widehat{\overline{\mathcal{H}}%
}(t)\{d\xi(t)-d{\small \widehat{\xi}}(t)\}]\\
& =\mathbb{E}[{\textstyle\int_{0}^{T}}\mathbb{E}[\widehat{\overline{\mathcal{H}}}(t)|\mathcal{G}_t]\{d\xi(t)-d{\small \widehat{\xi}}(t)\}]\\
& \leq 0,
\end{align*}
{\small where the last inequality holds because of the maximum condition
\eqref{eq3.11}. We conclude that
\[
J(\widehat{u},\xi)-J(\widehat{u},\widehat{\xi})\leq0.
\]
The proof is complete. \hfill$\square$ \bigskip}

\subsection{A necessary maximum principle}

{\small \noindent Since the concavity condition is not always satisfied, it is
useful to have a necessary condition of optimality where this condition is not
required. \noindent Suppose that a control }${\small (\widehat{u},\widehat
{\xi})}${\small $\in\mathcal{A}\times\mathcal{K}$ is an optimal pair and that
$(v,\zeta)${\normalsize $\in\mathcal{A}\times\mathcal{K}.$ Define }%
$u^{\lambda}=u+\lambda v$ and $\xi^{\lambda}=\xi+\lambda\zeta$, for a non-zero
sufficiently small $\lambda.$ \noindent Assume that $(u^{\lambda},\xi
^{\lambda})\in\mathcal{A}\times\mathcal{K}.$ For each given $t\in\lbrack0,T],$
let $\eta=\eta(t)$ be a bounded $\mathcal{G}_{t}$-measurable random variable,
let $h\in\lbrack T-t,T]$ and define }

{\small
\begin{equation}
v(s):=\eta\mathbf{1}_{\left[  t,t+h\right]  }(s);s\in\left[  0,T\right]  .
\label{eq4.1}%
\end{equation}
Assume that the \emph{derivative process} $Y(t)$, defined by $Y(t):=\tfrac
{d}{d\lambda}X^{u^{\lambda},\xi}(t)|_{\lambda=0}$ exists. Then we see that%
\[
Y(t)=%
{\textstyle\int_{0}^{t}}
\{\tfrac{\partial b}{\partial x}(t,s)Y(s)+\tfrac{\partial b}{\partial
u}(t,s)v(s)\}ds+%
{\textstyle\int_{0}^{t}}
\{\tfrac{\partial\sigma}{\partial x}(t,s)Y(s)+\tfrac{\partial\sigma}{\partial
u}(t,s)v(s)\}dB(s),
\]
and hence%
\begin{align}
dY(t)  &  =[\tfrac{\partial b}{\partial x}(t,t)Y(t)+\tfrac{\partial
b}{\partial u}(t,t)v(t)+%
{\textstyle\int_{0}^{t}}
(\tfrac{\partial^{2}b}{\partial t\partial x}(t,s)Y(s)+\tfrac{\partial^{2}%
b}{\partial t\partial u}(t,s)v(s))ds+%
{\textstyle\int_{0}^{t}}
(\tfrac{\partial^{2}\sigma}{\partial t\partial x}(t,s)Y(s)\nonumber\\
&  +\tfrac{\partial^{2}\sigma}{\partial t\partial u}(t,s)v(s))dB(s)]dt+(\tfrac
{\partial\sigma}{\partial x}(t,t)Y(t)+\tfrac{\partial\sigma}{\partial
u}(t,t)v(t))dB(t). \label{1.14}%
\end{align}
Similarly, we define the derivative process $Z(t):=\tfrac{d}{d\lambda}%
X^{u,\xi^{\lambda}}(t)|_{\lambda=0}$, as follows%
\[
Z(t)=%
{\textstyle\int_{0}^{t}}
\tfrac{\partial b}{\partial x}(t,s)Z(s)ds+%
{\textstyle\int_{0}^{t}}
\tfrac{\partial\sigma}{\partial x}(t,s)Z(s)dB(s)+%
{\textstyle\int_{0}^{t}}
h(t,s)d\zeta(s),
\]
which is equivalent to%
\begin{align}
dZ(t)  &  =[\tfrac{\partial b}{\partial x}(t,t)Z(t)+%
{\textstyle\int_{0}^{t}}
\tfrac{\partial^{2}b}{\partial t\partial x}(t,s)Z(s)ds]dt+\tfrac
{\partial\sigma}{\partial x}(t,t)Z(t)dB(t)\nonumber\\
&  +%
{\textstyle\int_{0}^{t}}
\tfrac{\partial^{2}\sigma}{\partial t\partial x}(t,s)Z(s)dB(s)dt+h(t,t)d\zeta
(t)+%
{\textstyle\int_{0}^{t}}
\tfrac{\partial h}{\partial t}(t,s)d\zeta(s)dt. \label{eq3.17}%
\end{align}
We shall prove the following theorem: }

\begin{theorem}
[Necessary maximum principle]{\small \vskip 0.2cm }

\begin{enumerate}
\item {\small For fixed $\xi{\small \in\mathcal{K}},$ suppose that
}${\small \widehat{u}}${\small $\in$ $\mathcal{A}$ is such that, for all
$\beta$ as in \eqref{eq4.1},
\begin{equation}
\tfrac{d}{d\lambda}J(\widehat{u}+\lambda\beta,\xi)|_{\lambda=0}=0\label{eq4.5}%
\end{equation}
and the corresponding solution }${\small \widehat{X}}${\small $(t),(\widehat
{p}(t),\widehat{q}(t,t))$ of (\ref{sde}) and (\ref{p}) exists. Then,
\begin{equation}
\mathbb{E[}\tfrac{\partial\mathbb{H}}{\partial u}(t)|\mathcal{G}%
_{t}]_{u=\widehat{u}(t)}=0.\label{eq4.6}%
\end{equation}
}

\item {\small Conversely, if \eqref{eq4.6} holds, then \eqref{eq4.5} holds.
}

{\small
}

\item {\small Similarly, for fixed }${\small \widehat{u}}$$\in$ $\mathcal{A}%
${\small , suppose that }${\small \widehat{\xi}}${\small $\in\mathcal{K}$ is
optimal. Then the following variational inequalities hold:
\begin{align}
\mathbb{E}[  \widehat{f_{1}}(t)+\widehat{p}(t)h(t,t)+%
{\textstyle\int_{t}^{T}}
\widehat{p}(s)\tfrac{\partial h}{\partial s}(s,t)ds|\mathcal{G}_{t}]
\leq0, \label{eq3.21}
\end{align}
and }
{\small
\begin{align}
\mathbb{E}[  \widehat{f_{1}}(t)+\widehat{p}(t)h(t,t)+%
{\textstyle\int_{t}^{T}}
\widehat{p}(s)\tfrac{\partial h}{\partial s}(s,t)ds|\mathcal{G}_{t}]
d\widehat{\xi}(t)=0. \label{eq3.22}
\end{align}
}
\end{enumerate}
\end{theorem}

{\small \noindent{Proof.} \quad For simplicity of notation we drop the
"hat" notation in the following.\newline Points 1-2 are direct consequence of Theorem 4.4 in Agram
\textit{et al} \cite{AOY2}. We proceed to prove point 3. Since $\widehat{u}$  is fixed we
drop the hat from the notation. Set
\begin{equation}%
\begin{array}
[c]{l}%
\tfrac{d}{d\lambda}J(\xi^{\lambda})|_{\lambda=0}
=\mathbb{E[}%
{\textstyle\int_{0}^{T}}
\{\tfrac{\partial f_{0}}{\partial x}(t)Z(t)dt+%
{\textstyle\int_{0}^{T}}
\tfrac{\partial f_{1}}{\partial x}(t)Z(t)d\xi(t)+%
{\textstyle\int_{0}^{T}}
f_{1}(t)d\zeta(t)+\tfrac{\partial g}{\partial x}(T)Z(T)].
\end{array}
\label{1.15}%
\end{equation}
Applying the It\^o formula, we get
\[%
\begin{array}
[c]{l}%
\mathbb{E[}\tfrac{\partial g}{\partial x}(T)Z(T)]=\mathbb{E[}p(T)Z(T)]\\
=\mathbb{E[}%
{\textstyle\int_{0}^{T}}
p(t)\{\tfrac{\partial b}{\partial x}(t,t)Z(t)+%
{\textstyle\int_{0}^{t}}
\tfrac{\partial^{2}b}{\partial t\partial x}(t,s)Z(s)ds\}dt\\
+%
{\textstyle\int_{0}^{T}}
p(t)\frac{\partial\sigma}{\partial x}(t,t)Z(t)dB(t)+%
{\textstyle\int_{0}^{T}}
p(t)(%
{\textstyle\int_{0}^{t}}
\frac{\partial^{2}\sigma}{\partial t\partial x}(t,s)Z(s)dB(s))dt\\
+%
{\textstyle\int_{0}^{T}}
p(t)h(t,t)d\zeta(t)+%
{\textstyle\int_{0}^{T}}
p(t)(%
{\textstyle\int_{0}^{t}}
\frac{\partial h}{\partial t}(t,s)d\zeta(s))dt\\
-%
{\textstyle\int_{0}^{T}}
Z(t)\tfrac{\partial\mathcal{H}}{\partial x}(t)dt-%
{\textstyle\int_{0}^{T}}
Z(t)\tfrac{\partial\overline{\mathcal{H}}}{\partial x}(t)d\xi(t)-%
{\textstyle\int_{0}^{T}}
Z(t)(%
{\textstyle\int_{t}^{T}}
\frac{\partial q}{\partial t}(t,s)dB(s))dt\\
+%
{\textstyle\int_{0}^{T}}
Z(t)q(t,t)dB(t)+%
{\textstyle\int_{0}^{T}}
Z(t)\tfrac{\partial\sigma}{\partial x}(t,t)q(t,t)dt].
\end{array}
\]
Therefore, from (\ref{1})-(\ref{3}), we obtain
\[%
\begin{array}
[c]{l}%
\mathbb{E}\left[  p(T)Z(T)\right]  \\
=\mathbb{E}[%
{\textstyle\int_{0}^{T}}
Z(t)\{\tfrac{\partial b}{\partial x}(t,t)p(t)+%
{\textstyle\int_{t}^{T}}
\left(  \tfrac{\partial^{2}b}{\partial s\partial x}(s,t)p(s)+{\mathbb{E}%
(D_{t}p(s)|\mathcal{F}_{t})}\tfrac{\partial^{2}\sigma}{\partial s\partial
x}(s,t)\right)  ds\}dt\\
+%
{\textstyle\int_{0}^{T}}
p(t)h(t,t)d\zeta(t)+%
{\textstyle\int_{0}^{T}}
(%
{\textstyle\int_{t}^{T}}
p(s)\frac{\partial h}{\partial s}(s,t)ds)d\zeta(t)-%
{\textstyle\int_{0}^{T}}
Z(t)\tfrac{\partial\mathcal{H}}{\partial x}(t)dt\\
-%
{\textstyle\int_{0}^{T}}
Z(t)\tfrac{\partial\overline{\mathcal{H}}}{\partial x}(t)d\xi(t)+%
{\textstyle\int_{0}^{T}}
Z(t)\tfrac{\partial\sigma}{\partial x}(t,t)q(t,t)dt].
\end{array}
\]
Using the definition of $\mathcal{H}$ and }$\overline{\mathcal{H}}$ {\small in
$(\ref{h})-(\ref{h1})$},%
\begin{equation}%
\begin{array}
[c]{l}%
\tfrac{d}{d\lambda}J(u,\xi^{\lambda})|_{\lambda=0}\\
=\mathbb{E[}%
{\textstyle\int_{0}^{T}}
\big\{p(t)h(t,t)+f_{1}(t)+%
{\textstyle\int_{0}^{T}}
{\textstyle\int_{t}^{T}}
p(s)\frac{\partial h}{\partial s}(s,t)ds\big\}d\zeta(t)].
\end{array}
\label{n}%
\end{equation}
Thus,
\begin{align*}
0 &  \geq\tfrac{d}{d\lambda}J(u,\xi^{\lambda})|_{\lambda=0}\\
&  =\mathbb{E[}%
{\textstyle\int_{0}^{T}}
\big\{p(t)h(t,t)+f_{1}(t)+{\textstyle\int_{t}^{T}}p(s)\tfrac{\partial
h}{\partial s}(s,t)ds\big \}d\zeta(t)],
\end{align*}
for all $\zeta$$\in\mathcal{K(\widehat{\xi})}$. \newline If we choose $\zeta$
to be a pure jump process of the form $\zeta(t)=\underset{0\leq t_{i}\leq
T}{\sum}\alpha(t_{i})$ where $\alpha(t_{i})>0$ is $\mathcal{G}_{t_{i}}%
$-measurable for all $t_{i}$, then $\zeta$$\in\mathcal{K(\widehat{\xi})}$ and
(\ref{n}) gives%
\[
\mathbb{E[(}f_{1}(t)+p(t)h(t,t)+%
{\textstyle\int_{t}^{T}}
p(s)\tfrac{\partial h}{\partial s}(s,t)ds)\alpha(t_{i})]\leq0\text{ for each
}t_{i}\text{ a.s.}%
\]
Since this holds for all such $\zeta$ with arbitrary $t_{i}$, we conclude that%
\[
\mathbb{E[(}f_{1}(t)+p(t)h(t,t)+%
{\textstyle\int_{t}^{T}}
p(s)\tfrac{\partial h}{\partial s}(s,t)ds)|\mathcal{G}_{t}]\leq0\text{ for
each }t\in\lbrack0,T]\text{ a.s.}%
\]
Finally, applying (\ref{n}) to $\zeta_{1}=\widehat{\xi}\in\mathcal{K(\widehat
{\xi})}$ and to $\zeta_{2}=-\widehat{\xi}\in\mathcal{K(\widehat{\xi})}$, we
get for all $t\in\lbrack0,T]$%
\[
\mathbb{E[(}f_{1}(t)+p(t)h(t,t)+%
{\textstyle\int_{t}^{T}}
p(s)\tfrac{\partial h}{\partial s}(s,t)ds)|\mathcal{G}_{t}]d\widehat{\xi
}(t)=0\text{ for each }t\in\lbrack0,T]\text{ a.s.}%
\]
\fproof

\section{Application to optimal harvesting with memory}
\subsection{Optimal harvesting with density-dependent prices}

{\small \noindent Let {\small $X^{\xi}(t)=X(t)$ be a given population
density (or cash flow) process, modelled by the following stochastic Volterra
equation:
\begin{equation}
X(t)=x_{0}+%
{\textstyle\int_{0}^{t}}
b_{0}(t,s)X(s)ds+{%
{\textstyle\int_{0}^{t}}
}\sigma_{0}(s)X(s)dB(s)-%
{\textstyle\int_{0}^{t}}
h(t,s)d\xi(s),\label{eq4.1a}%
\end{equation}
or, in differential form, }

{\small
\begin{equation}
\left\{
\begin{array}
[c]{l}%
dX(t)=b_{0}(t,t)X(t)dt+\sigma_{0}(t)X(t)dB(t)-h(t,t)\xi(t)\\
\quad\quad\quad+[\int_{0}^{t}\frac{\partial b_{0}}{\partial t}(t,s)X(s)ds-%
{\textstyle\int_{0}^{t}}
\frac{\partial h}{\partial t}(t,s)d\xi(s)]dt,\quad t\geq0.\\
X(0)=x_{0}.
\end{array}
\right.  \label{eq4.2}%
\end{equation}
We see that the dynamics of $X(t)$ contains a history or memory term
represented by the $ds$-integral$.$\newline We assume that $b_{0}(t,s)$ and
$\sigma_{0}(s)$ are given deterministic functions of $t$, $s$, with values in
$\mathbb{R}$, and that $b_{0}(t,s), h(t,s)$ are continuously
differentiable with respect to $t$ for each $s$ and $h(t,s)>0$. For
simplicity we assume that these functions are bounded, and the initial value
$x_{0}\in%
\mathbb{R}
$. 

We want to solve the following maximisation problem: 

\begin{problem}
{\small Find }${\small \widehat{\xi}}${\small $\in\mathcal{K},$ such that
\begin{equation}
\sup_{\xi}J(\xi)=J(\widehat{\xi}), \label{eq6.4a}%
\end{equation}
where
\begin{equation}
J(\xi)=\mathbb{E[}\theta X(T)+%
{\textstyle\int_{0}^{T}}
X(t)d\xi(t)]. \label{eq5.18}%
\end{equation}
}
\end{problem}

{\small \noindent Here $\theta=\theta(\omega)$ is a given $\mathcal{F}_{T}%
$-measurable square integrable random variable.\newline In this case the
Hamiltonian $\mathbb{H}$ takes the form
\begin{equation}%
\begin{array}
[c]{c}%
\mathbb{H}(t,x,p,q)=[b_{0}(t,t)xp+\sigma_{0}(t)xq+%
{\textstyle\int_{t}^{T}}
\tfrac{\partial b_{0}}{\partial s}(s,t)xp(s)ds-%
{\textstyle\int_{t}^{T}}
\tfrac{\partial h}{\partial s}(s,t)p(s)d\xi(s)]dt\\
+[x-h(t,t)p]d\xi(t).
\end{array}
\label{eq5.16}%
\end{equation}
Note that $\mathbb{H}$ is not concave with respect to $x$, so the sufficient
maximum principle does not apply. However, we can use the necessary maximum
principle as follows: The adjoint equation takes the form
\[
\left\{
\begin{array}
[c]{l}%
dp(t)=-\Big[p(t)b_{0}(t,t)+\sigma_{0}(t)q(t,t)+{\int_{t}^{T}}\frac{\partial
b_{0}}{\partial s}(s,t)p(s)ds\Big]dt+d\xi(t)+q(t,t)dB(t)\\
p(T)=\theta,
\end{array}
\right.
\]
equivalently%
\begin{equation}
p(t)=\theta+%
{\textstyle\int_{t}^{T}}
\{b_{0}(t,s)p(s)+\sigma_{0}(s)q(t,s)\}ds+%
{\textstyle\int_{t}^{T}}
d\xi(s)-%
{\textstyle\int_{t}^{T}}
q(t,s)dB(s). 
\end{equation}
\noindent The variational inequalities for an optimal control 
${\small \widehat{\xi}}$ and the corresponding ${\small \widehat{p}}$ are:
\begin{align}
&  \widehat{X}(t)-h(t,t)\widehat{p}(t)-%
{\textstyle\int_{t}^{T}}
\tfrac{\partial h}{\partial s}(s,t)\widehat{p}(s)ds\leq0,
\label{eq4.7a}\\
&  \text{ and }\nonumber\\
&  \big\{\widehat{X}(t)-h(t,t)\widehat{p}(t)-%
{\textstyle\int_{t}^{T}}
\tfrac{\partial h}{\partial s}(s,t)\widehat{p}(s)ds\big\}d\widehat
{\xi}(t)=0. \label{eq4.8a}%
\end{align}
}

{\small We have proved: }

\begin{theorem}
{\small Suppose $\hat{\xi}$ is an optimal control for Problem 4.1, with
corresponding solution $\hat{X}$ of \eqref{eq4.1a}. Then \eqref{eq4.7a} and
\eqref{eq4.8a} hold, i.e.
\begin{align}
&  \gamma_{0}(t,t)\widehat{p}(t)+{\textstyle\int_{t}^{T}}\tfrac{\partial h}{\partial s}(s,t)\widehat{p}(s)ds\geq \widehat{X}(t) \quad
a.s.,\quad t\in\lbrack0,T]\label{eq4.9a}\\
&  \text{ and }\nonumber\\
&  \big\{\gamma_{0}(t,t)\widehat{p}(t)+{\textstyle\int_{t}^{T}}\tfrac
{\partial h}{\partial s}(s,t)\widehat{p}(s)ds-\widehat
{X}(t)\big\}d\widehat{\xi}(t)=0. \label{eq4.10a}%
\end{align}
}
\end{theorem}

\begin{remark}
{\small The above result states that }${\small \widehat{\xi}}${\small $(t)$
increases only when
\begin{equation}
\gamma_{0}(t,t)\widehat{p}(t)+%
{\textstyle\int_{t}^{T}}
\tfrac{\partial h}{\partial s}(s,t)\widehat{p}(s)ds-\log(\widehat
{X}(t))=0. \label{eq4.11a}%
\end{equation}
Combining this with \eqref{eq4.7a} we can conclude that the optimal control can
be associated to the solution of a system of reflected forward-backward SVIEs
with barrier given by \eqref{eq4.9a}. }
\end{remark}
}
In particular, if we choose $h=1$  the variational inequalities become
\small{
\begin{align}
&  \widehat{p}(t) \geq \widehat{X}(t) \text{ for all } t \quad
a.s.,\quad t\in\lbrack0,T]\label{eq4.12a}\\
&  \text{ and }\nonumber\\
& \big\{\widehat{p}(t)-\widehat
{X}(t)\big\}d\widehat{\xi}(t)=0. \label{eq4.13a}%
\end{align}
}
\begin{remark}
This is a coupled system ($\widehat{X}(t),\widehat{p}(t))$ consisting of the solution $X(t)$ of the singularly controlled forward SDE 
\begin{equation}
X(t)=x_{0}+%
{\textstyle\int_{0}^{t}}
b_{0}(t,s)X(s)ds+{%
{\textstyle\int_{0}^{t}}
}\sigma_{0}(s)X(s)dB(s)-%
{\textstyle\int_{0}^{t}}
d\xi(s),\label{eq4.14a}%
\end{equation}
and the backward reflected SDE 
\begin{equation}
p(t)=\theta+%
{\textstyle\int_{t}^{T}}
\{b_{0}(t,s)p(s)+\sigma_{0}(s)q(t,s)\}ds-%
{\textstyle\int_{t}^{T}}
q(t,s)dB(s)+%
{\textstyle\int_{t}^{T}}
d\xi(s). \label{adjoint}%
\end{equation}
with barrier $\widehat{X}(t)$ and solution $\widehat{p}(t)$, if we choose $\xi=\widehat{\xi}$.
The optimal control is the process $\xi(t)$ which makes \eqref{eq4.12a} - \eqref{adjoint}  satisfied.
To the best of our knowledge such a forward-backward singularly controlled system has not been studied before. 
This is an interesting topic for future research.
 \end{remark}
 \color{black}
 
\subsection{Optimal harvesting with density-independent prices}

{\small \noindent\noindent Consider again equation (\ref{eq4.14a}) but now with
performance functional%
\[
J(\xi)=\mathbb{E[}\theta X(T)+%
{\textstyle\int_{0}^{T}}
\rho(t)d\xi(t)],
\]
for some positive deterministic function $\rho$. We want to find an optimal
}${\small \widehat{\xi}}${\small $\in\mathcal{K}$, such that }%
\[
\underset{\xi}{{\small \sup}}{\small J(\xi)=J(\widehat{\xi}).}%
\]
{\small In this case the Hamiltonian $\mathbb{H}$ gets the form
\[%
\begin{array}
[c]{c}%
\mathbb{H}(t,x,p,q)=[b_{0}(t,t)xp+\sigma_{0}(t)xq+%
{\textstyle\int_{t}^{T}}
\tfrac{\partial b_{0}}{\partial s}(s,t)xp(s)ds-%
{\textstyle\int_{t}^{T}}
\tfrac{\partial h}{\partial s}(s,t)p(s)d\xi(s)]dt\\
+[\rho(t)-h(t,t) p]d\xi(t).
\end{array}
\]
Note that $\mathbb{H}(x)$ is concave in this case. Therefore we can apply the
sufficient maximum principle here. The adjoint equation gets the form }

{\small
\begin{equation}
\left\{
\begin{array}
[c]{ll}%
dp(t) & =-\Big[p(t)b_{0}(t,t)+\sigma_{0}(t)q(t,t)+{\int_{t}^{T}}\frac{\partial
b_{0}}{\partial s}(s,t)p(s)ds\Big]dt+q(t,t)dB(t),\\
p(T) & =\theta,
\end{array}
\right.  \label{BSDE}%
\end{equation}
equivalently%
\[
p(t)=\theta+%
{\textstyle\int_{t}^{T}}
\{b_{0}(t,s)p(s)+\sigma_{0}(s)q(t,s)\}ds-%
{\textstyle\int_{t}^{T}}
q(t,s)dB(s).
\]
A closed form expression for $p(t)$ is given in the Appendix (Theorem 5.1).\\

\bigskip\noindent In this case the variational inequalities for an optimal
control ${\small \widehat{\xi}}$ and the corresponding
${\small \widehat{p}}$ are:
\begin{align}
&  \gamma_{0}(t,t)\widehat{p}(t)+%
{\textstyle\int_{t}^{T}}
\tfrac{\partial h}{\partial s}(s,t)\widehat{p}(s)ds\geq\rho
(t)\label{eq4.11}\\
&  \text{ and }\nonumber\\
&  \{\gamma_{0}(t,t)\widehat{p}(t)+%
{\textstyle\int_{t}^{T}}
\tfrac{\partial h}{\partial s}(s,t)\widehat{p}(s)ds-\rho
(t)\}d\widehat{\xi}(t)=0. \label{eq4.12}%
\end{align}
We have proved: }

\begin{theorem}
{\small Suppose $\hat{\xi}$ with corresponding solution }${\small \widehat{p}%
}${\small $(t)$ of the BSVIE \eqref{BSDE} satisfies the equations
\eqref{eq4.11} - \eqref{eq4.12}. Then ${\small \widehat{\xi}}$ is an
optimal control for Problem 4.1. }
\end{theorem}

\begin{remark}
{\small Note that \eqref{eq4.11} - \eqref{eq4.12} constitute a sufficient
condition for optimality. We can for example get this equation satisfied by
choosing $(\widehat{p}(t),\widehat{\xi}(t))$ as the solution of the BSVIE
\eqref{BSDE} reflected downwards at the barrier given by
\begin{align}
\gamma_{0}(t,t)\widehat{p}(t)+
{\textstyle\int_{t}^{T}}
\tfrac{\partial h}{\partial s}(s,t)\widehat{p}(s)ds-\rho(t)=0.
\end{align}
}
\end{remark}

\section{Appendix}

\begin{theorem}
{\small Consider the following linear BSVIE with singular drift
\begin{equation}
p(t)=\theta+%
{\textstyle\int_{t}^{T}}
\{b_{0}(t,s)p(s)+\sigma_{0}(s)q(t,s)\}ds+%
{\textstyle\int_{t}^{T}}
\tfrac{1}{X(s)}d\xi(s)-%
{\textstyle\int_{t}^{T}}
q(t,s)dB(s).
\end{equation}
}

{\small The first component $p(t)$ of the solution $(p(t),q(t))$ can be written in closed formula as
follows%
\[
p(t)=\mathbb{E}[\{\theta+\theta%
{\textstyle\int_{t}^{T}}
\Psi(t,s)ds+%
{\textstyle\int_{t}^{T}}
{\textstyle\int_{t}^{T}}
\Psi(t,s)\tfrac{1}{X(r)}d\xi(r)ds\}K(T)|\mathcal{F}_{t}],
\]
where
\begin{equation}
\Psi(t,r):=\Sigma_{n=1}^{\infty}b_{0}^{n}(t,r) \label{psi}%
\end{equation}
and $K(T)$ is given by%
\[
K(T)=\exp(%
{\textstyle\int_{0}^{T}}
\sigma_{0}(s)dB(s)-\tfrac{1}{2}%
{\textstyle\int_{0}^{T}}
\sigma_{0}^{2}(s)ds).
\]
}
\end{theorem}

\noindent\noindent{Proof.} \quad\ The proof is an extension of Theorem
3.1 in Hu and \O ksendal \cite{ho} to BSVIE with singular drift. \noindent
Define the measure $Q$ by
\[
dQ=M(T)dP\text{ on }\mathcal{F}_{T},
\]
where $M(t)$ satisfies the equation%
\[
\left\{
\begin{array}
[c]{ll}%
dM(t) & =M(t)\sigma_{0}(t)dB(t),\quad t\in\lbrack0,T],\\
M(0) & =1,
\end{array}
\right.
\]
which has the solution
\[
M(t):=\exp(%
{\textstyle\int_{0}^{t}}
\sigma_{0}(s)dB(s)-\tfrac{1}{2}%
{\textstyle\int_{0}^{t}}
\sigma_{0}^{2}(s)ds),\quad t\in\lbrack0,T].
\]
Then under the measure $Q$ the process
\begin{equation}
B_{Q}(t):=B(t)-%
{\textstyle\int_{0}^{t}}
\sigma_{0}(s)ds,\quad t\in\lbrack0,T] \label{e.def_bq}%
\end{equation}
is a $Q$-Brownian motion. }

{\small \noindent\noindent For all $0\leq t\leq r\leq T,$ define
\[
b_{0}^{1}(t,r)=b_{0}(t,r)\,,\quad b_{0}^{2}(t,r)=%
{\textstyle\int_{t}^{r}}
b_{0}(t,s)b_{0}(s,r)ds,
\]
and inductively
\[
b_{0}^{n}(t,r)=%
{\textstyle\int_{t}^{r}}
b_{0}^{n-1}(t,s)b_{0}(s,r)ds\,,\quad n=3,4,\cdots\,.
\]
Note that if $|b_{0}(t,r)|\leq C$ (constant) for all $t,r$, then by induction
on $n\in%
\mathbb{N}
:|b_{0}^{n}(t,r)|\leq\tfrac{C^{n}T^{n}}{n!},$ for all $t,r,n$. Hence,
\[
\Psi(t,r):=\Sigma_{n=1}^{\infty}|b_{0}^{n}(t,r)|<\infty,
\]
for all $t,r$. By changing of measure, we can rewrite equation (\ref{adjoint})
as
\begin{equation}
\text{ }p(t)=\theta+%
{\textstyle\int_{t}^{T}}
b_{0}(t,s)p(s)ds+%
{\textstyle\int_{t}^{T}}
X^{-1}(s)d\xi(s)-%
{\textstyle\int_{t}^{T}}
q(t,s)dB_{Q}(s),\quad0\leq t\leq T, \label{e.1.2}%
\end{equation}
where the process $B_{Q}$ is defined by \eqref{e.def_bq}. Taking the
conditional $Q$-expectation on $\mathcal{F}_{t}$, we get
\begin{align}
p(t)  &  =\mathbb{E}_{Q}[\theta+%
{\textstyle\int_{t}^{T}}
b_{0}(t,s)p(s)ds+%
{\textstyle\int_{t}^{T}}
X^{-1}(s)d\xi(s)|\mathcal{F}_{t}]\nonumber\\
&  =\tilde{F}(t)+%
{\textstyle\int_{t}^{T}}
b_{0}(t,s)\mathbb{E}_{Q}[p(s)|\mathcal{F}_{t}]ds+\mathbb{E}_{Q}[%
{\textstyle\int_{t}^{T}}
X^{-1}(s)d\xi(s)|\mathcal{F}_{t}],\quad0\leq t\leq T, \label{1.1.5}%
\end{align}
where
\[
\tilde{F}(s)=\mathbb{E}_{Q}[\theta|\mathcal{F}_{s}\,].
\]
Fix $r\in\lbrack0,t]$. Taking the conditional $Q$-expectation on
$\mathcal{F}_{r}$ of \eqref{1.1.5}, we get
\[
{\mathbb{E}}_{Q}\left[  \hat{p}(t)|{\mathcal{F}}_{r}\right]  =\tilde{F}(r)+%
{\textstyle\int_{t}^{T}}
b_{0}(t,s)\mathbb{E}_{Q}[p(s)|\mathcal{F}_{r}]ds+\mathbb{E}_{Q}[%
{\textstyle\int_{t}^{T}}
X^{-1}(s)d\xi(s)|\mathcal{F}_{r}],\quad r\leq t\leq T\,.
\]
Put
\[
\tilde{p}(s)={\mathbb{E}}_{Q}\left[  p(s)|{\mathcal{F}}_{r}\right]  ,\quad
r\leq s\leq T\,.
\]
Then the above equation can be written as
\[
\tilde{p}(t)=\tilde{F}(r)+%
{\textstyle\int_{t}^{T}}
b_{0}(t,s)\tilde{p}(s)ds+\mathbb{E}_{Q}[%
{\textstyle\int_{t}^{T}}
X^{-1}(s)d\xi(s)|\mathcal{F}_{r}],\quad r\leq t\leq T\,.
\]
Substituting $\tilde{p}(s)=\tilde{F}(r)+\int_{s}^{T}b_{0}(s,\alpha)\tilde
{p}(\alpha)d\alpha+{\mathbb{E}}_{Q}[%
{\textstyle\int_{s}^{T}}
X^{-1}(\alpha)d\xi(\alpha)|\mathcal{F}_{r}]$ in the above equation, we obtain
\[%
\begin{array}
[c]{ll}%
\tilde{p}(t) & =\tilde{F}(r)+%
{\textstyle\int_{t}^{T}}
b_{0}(t,s)\{\tilde{F}(r)+%
{\textstyle\int_{s}^{T}}
b_{0}(s,\alpha)\tilde{p}(\alpha)d\alpha+\mathbb{E}_{Q}[%
{\textstyle\int_{s}^{T}}
X^{-1}(\alpha)d\xi(\alpha)|\mathcal{F}_{r}]\}ds\\
& =\tilde{F}(r)+%
{\textstyle\int_{t}^{T}}
b_{0}(t,s)\tilde{F}(r)ds+%
{\textstyle\int_{t}^{T}}
b_{0}(t,s)\mathbb{E}_{Q}[%
{\textstyle\int_{s}^{T}}
X^{-1}(\alpha)d\xi(\alpha)|\mathcal{F}_{r}]ds\\
& +%
{\textstyle\int_{t}^{T}}
b_{0}^{(2)}(t,\alpha)\tilde{p}(\alpha)d\alpha,\quad r\leq t\leq T\,.
\end{array}
\]
Repeating this, we get by induction
\begin{align*}
\tilde{p}(t)  &  =\tilde{F}(r)+%
{\textstyle\sum_{n=1}^{\infty}}
{\textstyle\int_{t}^{T}}
b_{0}^{n}(t,\alpha)\tilde{F}(r)d\alpha+%
{\textstyle\sum_{n=1}^{\infty}}
{\textstyle\int_{t}^{T}}
b_{0}^{n}(t,\alpha)\mathbb{E}_{Q}[%
{\textstyle\int_{s}^{T}}
X^{-1}(\alpha)d\xi(\alpha)|\mathcal{F}_{r}]d\alpha\\
&  =\tilde{F}(r)+%
{\textstyle\int_{t}^{T}}
\Psi(t,\alpha)\tilde{F}(r)d\alpha+%
{\textstyle\int_{t}^{T}}
\Psi(t,\alpha)\mathbb{E}_{Q}[%
{\textstyle\int_{s}^{T}}
X^{-1}(\alpha)d\xi(\alpha)|\mathcal{F}_{r}]d\alpha.
\end{align*}
where $\Psi$ is defined by (\ref{psi}). Now substituting $\tilde{p}(s)$ in
(\ref{1.1.5}), for $r=t,$ we obtain
\begin{align*}
p(t)  &  =\tilde{F}(t)+%
{\textstyle\int_{t}^{T}}
\Psi(t,s)\tilde{F}(t)ds+%
{\textstyle\int_{t}^{T}}
\Psi(t,s)\mathbb{E}_{Q}[%
{\textstyle\int_{s}^{T}}
X^{-1}(\alpha)d\xi(\alpha)|\mathcal{F}_{t}]ds\\
&  =\mathbb{E}_{Q}[\theta+\theta%
{\textstyle\int_{t}^{T}}
\Psi(t,s)ds+%
{\textstyle\int_{t}^{T}}
\Psi(t,s)%
{\textstyle\int_{s}^{T}}
X^{-1}(\alpha)d\xi(\alpha)ds|\mathcal{F}_{t}]\\
&  =\mathbb{E}_{Q}[\theta+\theta%
{\textstyle\int_{t}^{T}}
\Psi(t,s)ds+%
{\textstyle\int_{t}^{T}}
\Psi(t,s)ds%
{\textstyle\int_{t}^{T}}
X^{-1}(\alpha)d\xi(\alpha)|\mathcal{F}_{t}].
\end{align*}
}

\thanks{{\bf Acknowledgements}. Nacira Agram and Bernt \O ksendal gratefully acknowledge the financial support provided by the Swedish Research Council grant (2020-04697) and the Norwegian Research Council grant (250768/F20), respectively.}

\end{document}